\documentclass[11pt]{amsart}

\usepackage{url}
\usepackage[vmargin=2.5cm, hmargin=2.5cm]{geometry}
\usepackage{fourier}
\usepackage{tikz}
\usepackage{amssymb}

\title{An overview of the computational aspects of nonunique factorization invariants}

\newtheorem{theorem}{Theorem}

\theoremstyle{definition}
\newtheorem{example}{Example}

\newtheorem{remark}[theorem]{Remark}

\author{P. A. Garc\'{\i}a-S\'anchez}
\address{Departamento de \'Algebra and CITIC-UGR, Universidad de Granada, E-18071 Granada, Espa\~{n}a}
\email{pedro@ugr.es}
\thanks{The author is supported by the projects MTM2010-15595, FQM-343,  FQM-5849 and FEDER funds. Thanks to Alfred Gerlondiger for his comments and suggestions, and for encouraging me to write this overview. Also thanks to Alfredo S\'anchez-R.-Navarro for his comments.}

\makeatletter
\def\@setauthors{%
  \begingroup
  \def\thanks{\protect\thanks@warning}%
  \trivlist
  \centering\footnotesize \@topsep30\p@\relax
  \advance\@topsep by -\baselineskip
  \item\relax
  \author@andify\authors
  \def\\{\protect\linebreak}%
  \authors%
  \ifx\@empty\contribs
  \else
    ,\penalty-3 \space \@setcontribs
    \@closetoccontribs
  \fi
  \endtrivlist
  \endgroup
}
\def\@settitle{\begin{center}%
  \baselineskip14\p@\relax
    \bfseries
  \@title
  \end{center}%
}
\makeatother

\begin{document}

\begin{abstract}
We give an overview of the existing algorithms to compute nonunique factorization invariants in finitely generated monoids.
\end{abstract}

\maketitle

\section{Introduction}

In this manuscript we give a general overview of the existing procedures to deal with nonunique factorization invariants. These methods have gained importance since they provide batteries of examples that can be used to understand how to prove theoretical results (or disprove ideas that we initially thought would hold). The algorithms are fed from the theory and in many cases from advances in integer linear programming. Thus in a sense, this is a wheel: theory produces algorithms that can be used to test new ideas, and these yield new results.

A \emph{semigroup} is a set with a binary associative operation. If a semigroup has an identity element (an element that operated with any other, in both sides, keeps the element unchanged), then we say that the semigroup is a \emph{monoid}. Let $(M,\cdot)$ be a monoid. An element $m\in M$ is a \emph{unit} if there exists $m'\in M$ such that $m \cdot m'=e=m'\cdot m$, where $e$ is the identity element of $M$. A monoid is \emph{reduced} if the only unit is the identity element. We are concerned with factorizations up to units, so we can at the very beginning remove the units from our monoid and consider that it is reduced. 

A monoid $M$ is \emph{commutative} if $m\cdot m'=m'\cdot m$ for all $m,m'\in M$. All monoids in this paper are commutative, and thus we will adopt additive notation in the following, and will use $0$ to the identity element.

A monoid $M$ is \emph{cancellative} if whenever $m+m'=m+m''$ for some $m, m',m''\in M$, we have $m'=m''$. If $(R,+\cdot)$ is a domain, then the underlying monoid $(R,\cdot)$ is commutative and cancellative. As with commutativity, we will also assume that our monoids are cancellative. 

Thus in what follows a monoid $M$ is meant to be commutative, cancellative and reduced. We denote $M^*=M\setminus \{0\}$.

An element $m$ in $M^*$ is said to be an \emph{atom} or \emph{irreducible} if whenever $m=m'+m''$ for some $m',m''\in M$, then either  $m'=0$ or $m''=0$ (recall that we are assuming that $M$ is reduced). Let $\mathcal A(M)$ denote the set of atoms of $M$. We say that $M$ is \emph{atomic} if every element $m\in M$ can be expressed as a sum of finitely many atoms. 

For a given set $X$, let $\mathcal F(X)$ be the free monoid on $X$, that is, the expressions of the form $\sum_{x\in X} \lambda_x x$ with $\lambda_x\in \mathbb N$ ($\mathbb N$ denotes the set of nonnegative integers), and all but finitely many $\lambda_x$ are zero. For $M$ an atomic monoid, denote by $\mathsf Z(M)=\mathcal F(\mathcal A(M))$. There is a natural monoid epimorphism 
\[\varphi: \mathsf Z(M)\to M,\ \varphi\Big(\sum_{a\in \mathcal A(M)}\lambda_a a\Big)= \sum_{a\in \mathcal A(M)}\lambda_a a.\]
Observe that many expressions of the form $\sum_{a\in \mathcal A(M)}\lambda_a a$ may correspond to the same element in $M$.  For $m\in M$, we define $\mathsf Z(m)=\varphi^{-1}(m)$. Every element in $\mathsf Z(m)$ is a \emph{factorization} of $m$. For $N\subseteq M$, we will write  $\mathsf Z(N)=\bigcup_{m\in N} \mathsf Z(m)$.

It may happen that the cardinality of $\mathsf Z(m)$ is one for all $m$ (and consequently $\varphi$ is an isomorphism and $M$ is a free monoid); in this case $M$ is said to be a \emph{factorial monoid}. It also may happen that there are finitely many factorizations for every element in the monoid $M$, and then we say that $M$ is a \emph{FF-monoid}. The \emph{length} of a factorization $\sum_{a\in \mathcal A(M)}\lambda_a a$ is $\sum_{a\in \mathcal A(M)}\lambda_a$. If for every element $m\in M$, all the lengths of its factorizations coincide, then we say that $M$ is a \emph{half-factorial} monoid; and if the set of possible lengths of factorizations are finite for every element, the the monoid is a \emph{BF-monoid} (see \cite{g-hk} for more details and properties of these monoids).

Observe that for computational aspects it is desirable that $M$ can be described in a ``finite'' way, and this happens in the case $M$ is an atomic monoid with finitely many atoms. In this setting, if the cardinality of $\mathcal A(M)$ is $e$, we can identify $\mathsf Z(M)$ with $\mathbb N^e$. As we are assuming $M$ is cancellative and reduced, this implies, that any two factorizations are incomparable with respect to the usual partial ordering in $\mathbb N^e$.  Dickson's lemma implies that $\mathsf Z(m)$ will have finitely many elements for any $m\in M$.

Transfer homomorphisms allow to study the arithmetical invariants (such as sets of lengths and catenary degree) of Krull and weakly Krull monoids in associated auxiliary monoids. In many cases these auxiliary monoids are finitely generated (see \cite{g-hk}). So, in these cases we will have FF-monoids, and we will be able to determine some properties using a computer.

Notice also that if we are assuming that $M$ is finitely generated, then according to \cite[Proposition 3.1]{fg}, we can assume that $M$ ``lives'' in $\mathbb Z^k\times \mathbb Z_{d_1}\times \cdots \times \mathbb  Z_{d_r}$. If $A=\{m_1,\ldots,m_e\}$ is the set of atoms of $M$, then $M=\langle A\rangle = \big\{ \sum_{i=1}^e n_i m_i\mid n_1,\ldots, n_t\in \mathbb N\big\}$. For $m\in M$ the set of factorizations of $m$ corresponds with the set of nonnegative integer solutions of the system of equations 
\[ (m_1 \mid \cdots \mid m_e) (x_1\ldots x_e)^T = m,\]
where the $m_i$'s are written in columns, and the last $r$ equations are in congruences modulo $d_1,\ldots, d_r$, respectively. In order to deal with these equations in congruences we can introduce auxiliary variables and then project to the original ones (see for instance \cite[Chapter 7]{fg}). The software \texttt{Normaliz} (\cite{normaliz}) can handle this kind of systems of equations. 

By removing equations in congruences, we then have a monoid that is, \emph{torsion free}, that is whenever $km =km'$ for $k$ a positive integer and $m, m'\in M$, we have $m=m'$. Every finitely generated commutative, cancellative, reduced and torsion free monoid is isomorphic to a submonoid of $\mathbb N^k$ for some positive integer $k$ (this is known in the literature as Grillet's Theorem, see for instance \cite[Theorem 3.11]{fg}).  A monoid with all these conditions is called an \emph{affine semigroup}. The set of atoms of an affine semigroup $M$ is $M^*\setminus (M^*+M^*)$, and it is the unique minimal generating system of $M$. So here minimal generators correspond with atoms (irreducibles).

We will give the arithmetic invariants in the scope of affine semigroups. This does not mean that the some of the methods reviewed can be used in a more general scope (even in an noncomputatonal framework), see for instance \cite{principal-ideal-noeth, phi-1, phi-2, phi-3}.

If $z$ and $z'$ are two factorizations of $m\in M$, then the pair $(z,z')$ is in the kernel of the monoid morphism $\varphi$, which in the setting of affine semigroups with atoms $\{m_1,\ldots, m_e\}$ can be written as
\[
\varphi: \mathbb N^e\to M,\ \varphi(n_1,\ldots, n_e)=n_1m_1+\cdots+ n_em_e.
\]
A \emph{presentation} $\sigma$ of $M$ is a generating system of $\ker \varphi=\{(x,y)\in \mathbb N^e\times \mathbb N^e\mid \varphi(x)=\varphi(y)\}$, that is, $\ker\varphi$ is the minimal congruence containing  $\sigma$.

\begin{remark}\label{remark-pres}
Notice that from the definition of presentation, if $\sigma$ is a presentation for $M$ and $z,z'$ are two factorizations of $m\in M$, then there exists a chain of factorizations $z_1,\ldots, z_r$ of $m$ such that 
\begin{itemize}
\item $z_1=z$, $z_r=z'$,
\item for every $i\in\{1,\ldots,r-1\}$ there exists $a_i,b_i,c_i\in\mathbb N^e$ such that $(z_i,z_{i+1})=(a_i+c_i,b_i+c_i)$ with either $(a_i,b_i)\in \sigma$ or $(b_i,a_i)\in\sigma$.
\end{itemize}
This idea actually catches the fact that $\ker\varphi$ is the least congruence containing $\sigma$, or in other words, it is the reflexive-symmetric-transitive closure of $\sigma$ compatible with addition.
\end{remark}

Hence knowing a presentation of $M$ (a generating set of $\ker \varphi$) allows us to know how to move from $z$ to $z'$, and consequently it will be a fundamental tool in the study factorizations of elements in affine semigroups. This is the case of tame and $\omega$-primality.

Recently it has been shown that some invariants are related to the calculation of the set of factorizations of a principal ideal, and if the monoid is full, there are specific procedures that significantly speed the process.

For numerical semigroups there are particularizations of the procedures based mainly in Ap\'ery sets, which avoid the use of linear integer programming, and work well for small generators. We will describe them when applicable.

This manuscript is meant to give a state of art of the implementations existing for the calculation of nonunique factorization invariants. We will simply explain the theory that supports these procedures, but will not describe deeply the functions used. We have implemented everything that is described here in the \texttt{GAP} (\cite{gap}) package \texttt{numericalsgps} (\cite{numericalsgps}), and thus it is part of this package (see the manual of the package for a description of the functions, examples and mode of operation). The reader interested in a full description and implementation of the algorithms can have a look at the source code available either in the \texttt{GAP} web page, or for the development version in \url{https://bitbucket.org/gap-system/numericalsgps} (the files containing the functions described here for numerical semigroups are in \texttt{ca\-te\-nary-tame.gi} and \texttt{contributions.gi};  those for affine semigroups are in \texttt{affine.gi}, both in the folder \texttt{gap}). The package tests availability of other packages (\cite{4ti2gap, 4ti2Interface, singular-gap, normalizinterface, gradedmodules}) that interact with \texttt{4ti2} (\cite{4ti2}), \texttt{Normaliz} (\cite{normaliz}) and \texttt{Singular} (\cite{Singular}). Depending on this availability, the package will use an specific method for the calculations. So we had to implement in some cases up to four methods for the same invariant depending on the extra software used (this is why there are several files with prefix \texttt{affine-extra} in the \texttt{gap} folder). 

\section{Presentations}

R\'edei proved in \cite{redei} that every finitely generated commutative monoid is finitely presented. In our setting this means that every affine semigroup admits a presentation with finitely many elements. Since then, many alternative and shorter proves have been published. We recall here one of these approaches.

Let $t$ be an unknown and $\mathbb k$ be a field. For $M$ an affine semigroup, define the \emph{semigroup ring} $\mathbb k[M]=\bigoplus_{m\in M} \mathbb k t^m$, where addition is performed component-wise and multiplication follows the rule $t^mt^{m'}=t^{m+m'}$.

Assume that $\{m_1,\ldots,m_e\}$ is a generating system of $M$. Herzog in \cite{herzog} proves that $\sigma$ is a presentation of $M$ if and only if the ideal $I_M=( X^a-X^b\mid (a,b)\in \sigma)$, where $I_M$ is the kernel of the ring homomorphism induced by 
\[
\mathbb k [x_1,\ldots,x_t] \to \mathbb k[M],\ x_i\mapsto t^{m_i}.
\]
Observe that for $n=(n_1,\ldots, n_k)$, we can write $t^m$ as $t_1^{n_1}\cdots t_k^{n_k}$ and in this way we can see $\mathbb k[M]$ as a subring of $\mathbb k[t_1,\ldots, t_k]$. In particular, we can compute a presentation of $M$ by using elimination: we start with the ideal $\big(x_1-t^{m_1},\ldots, x_e-t^{m_e}\big)\subseteq \mathbb k[x_1,\ldots, x_e,t_1,\ldots, t_k]$, and then eliminate the variables $t_1,\ldots, t_k$ to obtain $I_M$.

\begin{example}\label{pres}
Let us compute a presentation of $M=\langle (2,0),(0,2),(1,1),(2,1)\rangle$ with \texttt{singular}, \cite{Singular}. 

\begin{verbatim}
> ring r=0,(x,y,z,t,u,v),lp;
> ideal i = (x-u^2,y-v^2,z-u*v,t-u*v^2);
> eliminate(i,u*v);
_[1]=yz2-t2
_[2]=xt2-z4
_[3]=xy-z2
\end{verbatim}
This means that $I_M=\left(yz^2-t^2,x t^2-z^4, xy -z^2\right)$, and in light of Herzog's correspondence, the set  \[\big\{((0,1,2,0),(0,0,0,2)),((1,0,0,2),(0,0,4,0)),((1,1,0,0),(0,0,2,0))\big\}\] is a presentation for $M$.
\end{example}

A \emph{minimal presentation} of $M$ is a presentation that cannot be refined to another presentation of $M$, that is, it is minimal with respect to set inclusion (it turns out that it is also minimal with respect to cardinality; see \cite[Corollary 9.5]{fg}). 

\begin{example}\label{pres-min}
The presentation in Example \ref{pres} is not minimal. If we want to obtain a minimal presentation with \texttt{singular} additional work is needed.
\begin{verbatim}
> ring r=0,(x,y,z,t,u,v),(wp(2,2,2,3),lp(2));
// ** redefining r **
> ideal i = (x-u**2,y-v**2,z-u*v,t-u*v**2);
> ideal j=eliminate(i,u*v);
> minbase(j);
_[1]=xy-z2
_[2]=yz2-t2
\end{verbatim}
\end{example}

Given $m\in M$, we define $\nabla_m$ as the graph with vertices $\mathsf Z(m)$ and $zz'$ is an edge if $z\cdot z'\neq 0$ (dot product). An element $m$ is a \emph{Betti element} of $M$ if the graph $\nabla_m$ is not connected. We will denote by $\mathrm{Betti}(M)$ the set of Betti elements of $M$.

The sets of vertices of the connected components of $\nabla_m$ are also known as $\mathcal R$-classes of $\mathsf Z(m)$. The following method can be used to produce all minimal presentations (up to arrangement of the pairs and symmetry) of $M$; see for instance \cite[Chapter 9]{fg}. 

\begin{itemize}
\item For all $m\in M$, if $\nabla_m$ is connected, then set $\sigma_m=\emptyset$. If not, let $R_1,\ldots, R_r$ be the different $\mathcal R$-classes of $\mathsf Z(m)$. Consider any tree $T$ with vertices $R_1,\ldots, R_r$. For each $i\in\{1,\ldots,r\}$ take $r_i\in R_i$. Set $\sigma_m= \{ (z_i,z_j)\mid R_iR_j \hbox{ is an edge of } T\}$ (for instance, one might take $\sigma_m=\{(z_1,z_2),(z_1,z_3),\ldots, (z_1,z_r)\}$).

\item The set $\sigma=\bigcup_{m\in M} \sigma_m$ is a minimal presentation of $M$.   
\end{itemize}

It follows that the set of Betti elements of $M$ has finite cardinality. And that the cardinality of a (any) minimal presentation is $\sum_{b\in\mathrm{Betti}(M)} (\mathrm{ncc}(\nabla_b)-1)$, where $\mathrm{ncc}(\nabla_b)$ stands for the number of connected components of $\nabla_b$. This formula holds for every atomic monoid having the ascending chain on principal ideals (\cite[Corollary 1]{principal-ideal-noeth}).

\begin{example}
Let $M$ be as in Example \ref{pres}. Since any presentation contains a minimal presentation, we have that $\mathrm{Betti}(M)\subseteq \{(2,4),(2,2),(4,4)\}$. We use the \texttt{GAP} (\cite{gap}) package \texttt{numericalsgps} (\cite{numericalsgps}) to calculate the $\mathcal R$-classes of each of these elements.
\begin{verbatim}
gap> RClassesOfSetOfFactorizations(
 FactorizationsVectorWRTList([4,4],[[2,0],[0,2],[1,1],[1,2]]));
[ [ [ 0, 0, 4, 0 ], [ 1, 0, 0, 2 ], [ 1, 1, 2, 0 ], [ 2, 2, 0, 0 ] ] ]
gap> RClassesOfSetOfFactorizations(
 FactorizationsVectorWRTList([2,4],[[2,0],[0,2],[1,1],[1,2]]));
[ [ [ 0, 1, 2, 0 ], [ 1, 2, 0, 0 ] ], [ [ 0, 0, 0, 2 ] ] ]
gap> RClassesOfSetOfFactorizations(
 FactorizationsVectorWRTList([2,2],[[2,0],[0,2],[1,1],[1,2]]));                      
[ [ [ 1, 1, 0, 0 ] ], [ [ 0, 0, 2, 0 ] ] ]
\end{verbatim}
It follows that $\mathrm{Betti}(M)=\{(2,2),(2,4)\}$ (this also follows from Example \ref{pres-min}).

The function \texttt{FactorizationsVectorWRTList} either uses \cite{cd}, or if available \cite{normaliz} or \cite{4ti2} through the packages \texttt{NormalizInterface} (\cite{normalizinterface}) or either \texttt{4ti2gap} (\cite{4ti2gap}) or \texttt{4ti2Interface} (\cite{4ti2Interface}).
\end{example}

We will see that the catenary degree and the Delta sets are ``ruled'' by a minimal presentation.

\section{Ap\'ery sets}\label{sec:apery}

Let $M$ be an affine semigroup generated by $\{m_1,\ldots,m_e\}$. Let $m\in M$. The \emph{Ap\'ery set} of $m$ in $M$ is the set 
\[
\mathrm{Ap}(M,m)=\{ m'\in M\mid m'-m\not\in M\}.
\]
Ap\'ery sets can be defined in a more general setting. If our monoid fulfills the ascending chain condition on principal ideals, then every for every $m'\in M$ there exits unique $(w,k)\in \mathrm{Ap}(M,m)\times \mathbb N$ such that $m'=km+w$ (see \cite{principal-ideal-noeth}).

If $M$ is a numerical semigroup, then the cardinality of $\mathrm{Ap}(M,m)$ has exactly $m$ elements. Moreover, if $b\in \mathrm{Betti}(M)$, then $b=m_i+w$ with $i\in\{2,\ldots,e\}$ and $w\in \mathrm{Ap}(M,m_1)\setminus\{0\}$ (see for instance \cite[Proposition 8.19]{ns}).  

As minimal presentations are crucial for studying factorizations, the above paragraph implies that Ap\'ery sets are also important in our study particularized to the numerical semigroup setting.

\section{Graver bases}

Let $M$ be an affine semigroup, $M\subseteq \mathbb N^k$ generated by $\{m_1,\ldots, m_e\}$.

We have seen that a minimal presentation is a minimal generating system of $\ker\varphi$ as a congruence. It turns out that $\ker\varphi$ is not only a congruence, but an affine semigroup itself, and thus it admits a unique minimal generating system, which we denote by $\mathcal I(M)$. It follows easily that $\mathcal I(M)$ corresponds with the pairs $(x,y)=((x_1,\ldots, x_e),(y_1,\ldots,y_e))\in \mathbb N^e\times \mathbb N^e\setminus\{(0,0)\}$ that are minimal (with respect to the usual product order) solutions of 
\[
(m_1|\cdots|m_e|-m_1|\cdots|-m_e)(x\mid y)^T=0,
\]
because if $(x,y)\in \ker\varphi$, then $x_1m_1+\cdots +x_em_e= y_1m_1+\cdots +y_em_e$. Moreover, there exists $(x_1,y_1),\ldots, (x_s,y_s)\in \mathcal I(M)$ such that $(x,y)=(x_1,y_1)+\cdots + (x_s,y_s)$. That is, every pair of factorizations of the same element can be expressed as a sum of factorizations of some specific elements. Indeed, we will say that $m\in M$ is \emph{primitive} if there exists $x,y\in \mathsf Z(m)$ such that $(x,y)\in \mathcal I(M)$.

In particular, $\mathcal I(M)$ is a presentation of $M$, though in general with a lot of redundancy.

Notice that if $\mathbf e_i$ is the $i$th row of the identity $e\times e$ matrix, then $(\mathbf e_i,\mathbf e_i)\in \mathcal I(M)$ for all $i\in\{1,\ldots, e\}$. 

On $\mathbb Z^e$ define the order $(x_1,\ldots, x_e)\sqsubseteq (y_1,\ldots, y_e)$ if for all $i\in\{1,\ldots,e\}$, $x_iy_i\ge 0$ and $|x_i|\le |y_i|$. Also, for $x\in \mathbb Z^e$ set $x^+$ and $x^-$ to be the unique elements in $\mathbb N^e$ such that $x=x^+-x^-$ and $x^+\cdot x^-=0$. It turns out that $x\sqsubseteq y$ if and only if $(x^+,x^-)\le (y^+,y^-)$ (usual partial ordering). 

Let $H$ be a subgroup of $\mathbb Z^e$. A \emph{Graver basis} of $H$ is a set of minimal nonzero elements of $H$ with respect to $\sqsubseteq$. 

Notice that the set of integer solutions of 
\[
(m_1\mid \cdots \mid m_e)x^T=0
\]
defines a subgroup $H_M$ of $\mathbb Z^e$. In fact $(x,y)\in \ker\varphi$ if and only if $x-y\in H_M$ (this is a rephrasing of the necessity condition in \cite[Proposition 1.4]{fg}). From a Graver basis $G$ of $H_M$ we can easily compute 
\[\mathcal I(M)=\left\{(x^+,x^-)\mid x\in G\right\}\cup \big\{(\mathbf e_i, \mathbf e_i)\mid i\in\{1,\ldots, e\}\big\}.\]

\begin{example}
Let us go back to $M$ in Examples \ref{pres} and \ref{pres-min}.
\begin{verbatim}
gap> GraverBasis4ti2(["mat",TransposedMat([[2,0],[0,2],[1,1],[1,2]])]);
[ [ 1, 0, -4, 2 ], [ 0, 1, 2, -2 ], [ 1, 1, -2, 0 ], [ 1, 2, 0, -2 ] ]
\end{verbatim}
The output of \texttt{4ti2} does not print an element and its negation. Hence a Graver basis of $H_M$ consists in 8 elements and $\mathcal I(M)$ has 8+4 elements.
\end{example}

We will see that some nonunique factorization invariants depend on the factorizations of the primitive elements of $M$.

\section{Block monoids}

Let $G$ be an Abelian group. And let $g_1,\ldots, g_k\in G$. A \emph{zero-sum sequence} is an expression of the form $n_1g_1+\cdots +n_kg_k=0$. The \emph{length} of this sequence is $n_1+\cdots+n_k$.  We say that a zero sum sequence is \emph{minimal} if there is no other zero-sum sequence $n_1'g_1+\cdots +n_k'g_k=0$ such that $0\neq (n_1',\ldots, n_k')\lneq (n_1,\ldots, n_k)$. The set of zero-sum sequences is clearly a monoid, actually it can be identified as a submonoid of $\mathbb N^k$ and it is generated by the minimal zero-sum sequences (indeed it is a full affine semigroup). We will denote the set of zero-sum sequences in $g_1,\ldots, g_k$ by $\mathcal B(\{g_1,\ldots,g_k\})$. 

Since $G$ is an Abelian group, it is then isomorphic to $\mathbb Z_{d_1}\times \cdots \times \mathbb Z_{d_r} \times \mathbb Z^l$ for some $d_1,\ldots, d_r, l\in\mathbb N$. Hence we can identify the elements $g_1,\ldots,g_k$ with elements in  $\mathbb Z_{d_1}\times \cdots \times \mathbb Z_{d_r} \times \mathbb Z^l$.  Hence $\mathcal B(\{g_1,\ldots, g_k\})$ corresponds with the set of nonnegative integer solutions of the system of $r+l$ equations and $k$ unknowns
\[
(g_1\mid \cdots \mid g_k)x =0 \in \mathbb Z_{d_1}\times \cdots \times \mathbb Z_{d_r} \times \mathbb Z^l
\]
(the first $r$ equations are in congruences modulo $d_1,\ldots, d_r$, respectively). The set of solutions of this system of equations can be computed via \texttt{Normaliz} (\cite{normaliz}).

The \emph{Davenport constant} is the supremum (in this setting maximum) of the lengths of minimal zero-sum sequences.

\begin{example}
We can compute the block monoid associated to $\mathbb Z_2^2$ in the following way using \texttt{numericalsgps}.
\begin{verbatim}
gap> m2:=[[0,1],[1,0],[1,1]];;                                      
gap> a:=AffineSemigroup("equations",[TransposedMat(m2),[2,2]]);;
gap> GeneratorsOfAffineSemigroup(a);
[ [ 0, 0, 2 ], [ 0, 2, 0 ], [ 1, 1, 1 ], [ 2, 0, 0 ] ]
\end{verbatim}
Observe that we are omitting $(0,0)$ and that the second argument of \texttt{AffineSemigroup} is a matrix whose columns are the elements in $(\mathbb Z_2^2)^*$ and a list indicating the equations that are congruences with the respective modules. The Davenport constant in this case is 3.
\end{example}

Many factorization properties of monoids can be derived from the factorization properties (or bounded in some cases) of the block monoid of their class groups (see \cite{g-hk}). This is why these affine semigroups are relevant in the study of nonunique factorization invariants.

\section{Denumerant and maximal denumerant}

We have already mentioned that for an affine semigroup $M$ and $m\in M$, the set $\mathsf Z(m)$ has finitely many elements. The \emph{denumerant} of $m$ is precisely the cardinality of $\mathsf Z(m)$. There is a wide amount of literature devoted to the study of denumerants of elements in numerical semigroups, indeed few formulas are known, and just for some particular families of monoids (\cite{alfonsin} is a nice reference for the reader interested in this topic). 

Of course the bigger an integer in a numerical semigroup is, the larger is its denumerant, and thus it is not bounded. What is indeed astonishing is that the maximal denumerant is bounded for numerical semigroups. The \emph{maximal denumerant} of $m$ in $M$ is the number of elements in $\mathsf Z(m)$ with maximal length. If $M$ is a numerical semigroup, set the maximal denumerant of $M$ as the maximum of the maximal denumerants of elements of $M$. Bryant and Hamblin give  in \cite{max-den} a procedure to compute the maximal denumerant of any numerical semigroup.

\begin{example}
The semigroup $\langle 3,5,7\rangle$ has maximal denumerant 2. 
\begin{verbatim}
gap> s:=NumericalSemigroup(3,5,7);;
gap> MaximalDenumerantOfNumericalSemigroup(s);
2
gap> List(Intersection([0..100],s),                              
> x->Length(FactorizationsElementWRTNumericalSemigroup(x,s)));
[ 1, 1, 1, 1, 1, 1, 1, 2, 1, 2, 2, 2, 3, 2, 3, 3, 3, 4, 4, 4, 4, 5, 5, 5, 6, 
  6, 6, 7, 7, 7, 8, 8, 9, 9, 9, 10, 10, 11, 11, 12, 12, 12, 14, 13, 14, 15, 
  15, 16, 16, 17, 17, 18, 19, 19, 20, 20, 21, 22, 22, 23, 24, 24, 25, 26, 26, 
  27, 28, 29, 29, 30, 31, 31, 33, 33, 34, 35, 35, 37, 37, 38, 39, 40, 41, 41, 
  43, 43, 44, 46, 46, 47, 48, 49, 50, 51, 52, 53, 54, 55 ]
\end{verbatim}
\end{example}

\section{Length based invariants}

Let $M$ be an affine semigroup generated by $\{m_1,\ldots, m_e\}$. Take $m\in M$ and  $x=(x_1,\ldots,x_e)\in \mathsf Z(m)$. Recall that the length of $x$ is defined as 
\[
\lvert x \rvert = x_1+\dots+x_e.
\]  
The \emph{set of lengths of factorizations} of $m$ is 
\[
\mathsf L(m)=\big\{ \lvert x\rvert \mid x\in\mathsf Z(m)\big\}.
\]

Since $\mathsf Z(m)$ has finitely many elements, so has $\mathsf L(m)$. This means that affine semigroups are BF-monoids. 

Recall that a monoid is half factorial if the cardinality of $\mathsf L(s)$ is one for all $s\in S$. This concept was introduced for domains in \cite{zaks}. 

From Remark \ref{remark-pres} it easily follows that $M$ is half factorial if and only if for every $(a,b)$ in a minimal presentation of $M$ we have $|a|=|b|$ (see \cite{elasticity}). Thus we can determine whether or not an affine semigroup is half-factorial.

\begin{example}
In Example \ref{pres-min}, since $((1,2,0,0),(0,0,0,2))$ belongs to a minimal presentation of $M$, we deduce that $M$ is not half-factorial.
\end{example}

\subsection{Elasticity}
One of the first nonunique factorization invariants that appeared in the literature was the elasticity (introduced in \cite{valenza}). It was meant to measure how far is a monoid from being half factorial. 

Take $m$ in an affine semigroup $M$. The \emph{elasticity} of $m$, $\rho(m)$, is defined as 
\[
\rho(m)=\frac{\sup \mathsf L(m)}{\min \mathsf L(m)}.
\]
Since $\mathsf L(m)$ has finitely many elements, the supremum in the numerator is indeed a maximum. 
The elasticity of $M$ is defined as 
\[
\rho(M)=\sup\big\{ \rho(m) \mid {m\in M} \big\}.
\]
It is not hard to show (see \cite{elasticity}) that 
\[
\rho(M)=\max\left\{ \frac{|a|}{|b|} ~\Big|~ (a,b)\in \mathcal I(M)\right\}.
\]
Hence by computing a Graver basis of $H_M$ we can calculate the elasticity of $M$. Philipp in his thesis, and published later in \cite{phi-2}, did a great improvement in the computation of the elasticity: he showed that we only have to consider elements  $(a,b)\in \mathcal I(M)$ with $a\neq b$ and with minimal support (indices of nonzero coordinates). These elements are known in the literature as circuits, and we can use \cite[Lemma 8.8]{binomials} to calculate them by means of determinants.  The calculation of a Graver basis is in general a hard problem, while computing determinants is affordable. Thus Philipp theoretical contribution became a considerable speed up in combination with Eisendbud and Sturmfels method for computing circuits.
\[
\rho(M)=\max\left\{ \frac{|a|}{|b|} ~\Big|~ (a,b)\hbox{ circuit of }\ker\varphi\right\}.
\]

\begin{example}
Let us compute $\rho\left(\mathcal B\left(\mathbb Z_2^3\right)\right)$. 
\begin{verbatim}
gap> m:=[[0,0,1],[0,1,0],[0,1,1],[1,0,0],[1,0,1],[1,1,0],[1,1,1]];;
gap> a:=AffineSemigroup("equations",[TransposedMat(m),[2,2,2]]);;
gap> ElasticityOfAffineSemigroup(a);
2
\end{verbatim}

\end{example}

\subsection{Delta sets}
Another way to measure how far we are from half factoriality, is to determine how distant are the different lengths of factorizations. This is the motivation for the following definition.

Let as above $m$ be an element in the affine semigroup $M$. 
Assume that $\mathsf L(m)=\{l_1<\cdots < l_r\}$. Define the \emph{Delta set} of $m$ as 
\[
\Delta(s)=\{ l_2-l_1,\ldots, l_r-l_{r-1}\},
\]
and if $r=1$, $\Delta(m)=\emptyset$. The Delta set of $M$ is defined as 
\[
\Delta(M)=\bigcup_{m\in M} \Delta(m).
\]
So, the bigger $\Delta(M)$ is, the farther is $M$ from begin half factorial.

Recall that $(x,y)\in \ker\varphi$ if and only if $x-y\in H_M$. Indeed, it is not hard to show that $H_M$ is generated as a group by the differences of the pairs in a presentation of $M$.  From this, one can prove that 
\[
\min\Delta(M)=\gcd\Delta(M)
\]
(\cite[Proposition 1.4.4]{g-hk}).

By using the idea expressed in Remark \ref{remark-pres}, it can be shown that the maximum of the distances between lengths of factorizations is reached in a Betti element of $M$ (\cite[Theorem 2.5]{deltas}):
\[
\max\Delta(M)=\max\big\{ \max\Delta(b)\mid b\in \mathrm{Betti}(M)\big\}.
\]
This gives us an interval where the elements in $\Delta(M)$ must be, but it is far from being a procedure to compute the whole set $\Delta(M)$. 

For numerical semigroups, it is known that the sets of distances between consecutive lengths of factorizations are eventually periodic (\cite{delta-per}) and a bound for this periodicity is given. This bound was improved in \cite{delta-comp}. Hence we can compute the Delta sets of the elements up to this bound (a dynamic version of this procedure is resented in \cite{b-on-p}). The problem is that this bound can be huge.

\begin{verbatim}
gap> s:=NumericalSemigroup(701,902,1041); 
<Numerical semigroup with 3 generators>
gap> DeltaSetOfNumericalSemigroup(s);
[ 1, 2, 3, 4, 5, 6, 11, 17 ]
gap> DeltaSetPeriodicityBoundForNumericalSemigroup(s);
313436
\end{verbatim}

Recently in \cite{delta-tres} a procedure that runs as fast as Euclid's extended algorithm has been presented for numerical semigroups with embedding dimension three (and not symmetric, though the algorithm seems to work also for symmetric numerical semigroups).

O'Neil in \cite{hilbert} gives new tools based on Hilbert functions that probably will yield procedures for the computation of  $\Delta(M)$ for an arbitrary affine semigroup $M$.

\section{Distance based invariants}
Observe that length based invariants cannot describe the behavior of factorizations in half-factorial monoids. To measure how spread are the factorizations, we first need a distance. 

For $x=(x_1,\ldots, x_e),y=(y_1,\ldots, y_e)\in\mathbb N^e$, define the infimum of $x$ and $y$ as 
\[
x\wedge y=(\min\{x_1,y_1\},\ldots, \min\{x_p,y_p\})
\]
(if we think in additive notation and $x$ and $y$ are factorizations of an element, then $x\wedge y$ translates to greatest common divisor).

The \emph{distance} between $x$ and $y$ is defined as 
\[
\mathrm d(x,y)=\max\{ \lvert x-(x\wedge y)\rvert, \lvert y-(x\wedge y)\rvert\}
\]
(equivalently $\mathrm d(x,y)=\max\{ \lvert x\rvert, \lvert y\rvert\}-\lvert x\wedge y\rvert$).

\subsection{Catenary degree}

We start with an example that illustrates the idea of catenary degree.

\begin{example}\label{ejemplo-cat-estacas}
The factorizations of $66\in \langle 6,9,11\rangle$ are \[\mathsf Z(66)=\big\{ (0, 0, 6 ), ( 1, 3, 3 ), ( 2, 6, 0 ), (4, 1, 3 ), ( 5, 4, 0 ),( 8, 2, 0),( 11, 0, 0 ) \big\}.\] The distance between $(11,0,0)$ and $(0,0,6)$ is $11$. 

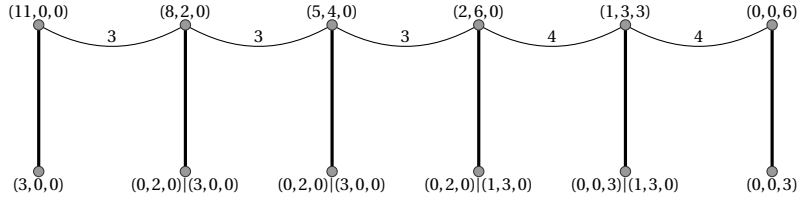
\begin{center}
\begin{tikzpicture}[y=.3cm, x=.3cm,font=\sffamily,, every  node/.style={scale=.65},scale=.65]
\draw[very 	thick, black] (0,0) -- (0,10);
\draw (0,10) to[bend right] (10,10)  to[bend right] (20,10) to[bend right] (30,10) to[bend right] (40,10) to[bend right] (50,10);
\draw[very thick, black] (10,10) -- (10,0);
\draw[very thick, black] (20,10) -- (20,0);
\draw[very thick, black] (30,10) -- (30,0);
\draw[very thick, black] (40,10) -- (40,0);
\draw[very thick, black] (50,10) -- (50,0);

\filldraw[fill=black!40,draw=black!80] (0,0) circle (3pt)    node[anchor=north] {$(3,0,0)$};

\filldraw[fill=black!40,draw=black!80] (0,10) circle (3pt)    node[anchor=south] 
{$( 11,0,0 )$};

\filldraw[fill=black!40,draw=black!80] (10,10) circle (3pt)    node[anchor=south] 
{$( 8,2,0 )$};

\filldraw[fill=black!40,draw=black!80] (10,0) circle (3pt)    node[anchor=north] 
{$(0,2,0) | (3,0,0) $};

\filldraw[fill=black!40,draw=black!80] (20,10) circle (3pt)    node[anchor=south] 
{$( 5,4,0 )$};

\filldraw[fill=black!40,draw=black!80] (20,0) circle (3pt)    node[anchor=north] 
{$(0,2,0) | (3,0,0) $};

\filldraw[fill=black!40,draw=black!80] (30,10) circle (3pt)    node[anchor=south] 
{$( 2,6,0 )$};

\filldraw[fill=black!40,draw=black!80] (30,0) circle (3pt)    node[anchor=north] 
{$(0,2,0) | (1,3,0) $};

\filldraw[fill=black!40,draw=black!80] (40,10) circle (3pt)    node[anchor=south] 
{$( 1,3,3 )$};

\filldraw[fill=black!40,draw=black!80] (40,0) circle (3pt)    node[anchor=north] 
{$(0,0,3) | (1,3,0) $};

\filldraw[fill=black!40,draw=black!80] (50,0) circle (3pt)    node[anchor=north] 
{$( 0,0,3 )$};

\filldraw[fill=black!40,draw=black!80] (50,10) circle (3pt)    node[anchor=south] {$( 0,0,6 )$};

\node [below] at (5,10) {$3$};

\node [below] at (15,10) {$3$};

\node [below] at (25,10) {$3$};

\node [below] at (35,10) {$4$};

\node [below] at (45,10) {$4$};

\end{tikzpicture}
\end{center}
In the above picture the factorizations are depicted in the top of a post, and they are linked by a ``catenary'' labeled with the distance between two consecutive sticks. On the bottom we have drawn the factorizations removing the 
common part with the one on the left and that of the right, respectively. So we have linked $(11,0,0)$ and $(0,0,6)$ with a chain of factorizations, and every two consecutive nodes in the chain are at most at distance 4. This is in fact the best we can do in this example. We are not caring about the length of the sequence, but about how closer are two consecutive elements in the chain.
\end{example}

Let $M$ be an affine semigroup, and take $m\in M$. Let $x,y\in \mathsf Z(m)$ and let $N$ be a nonnegative integer. An \emph{$N$-chain}\index{chain of factorizations} joining $x$ and $y$ is a sequence $x_1,\ldots, x_k\in\mathsf Z(m)$ such that 
\begin{itemize}
\item $x_1=x$, $x_k=y$,
\item for all $i\in\{1,\ldots, k-1\}$, $\mathrm d(x_i,x_{i+1})\le N$.
\end{itemize}

The \emph{catenary degree}\index{catenary degree} of $m$, denoted $\mathsf c(m)$, is the least $N$ such that for any two factorizations $x,y\in\mathsf Z(m)$, there is an $N$-chain joining them. The catenary degree of $M$, $\mathsf c(M)$, is defined as 
\[
\mathsf c(M)=\sup\{\mathsf c(m)\mid m\in M\}.
\]

The calculation of $\mathsf c(m)$ can be performed in the following way. We consider the complete graph with vertices the factorizations of $m$, and edges labeled with the distances between their ends. Then we pick an edge with the largest label, and if it is not a bridge, then we remove it. We keep doing so, until we arrive to a bridge. The label of this bridge is $\mathsf c(m)$.

\begin{example}
As an illustration of the above procedure, consider $77\in S = \langle 10,11, 23, 35 \rangle$.  In the following figure we see that we can remove the edge with label 6, meaning that in order to go from $(0,7,0,0)$ to $(2,1,2,0)$ we can first go to $(2,2,0,1)$ and then to $(2,1,2,0)$, and the distances in this walk between two consecutive nodes are less than 6. Then we remove the edge labeled with 5. But we cannot remove the edge joining $(1,4,1,0)$ and $(0,7,0,0)$ since it is a bridge (we can remove the other labeled with 3).

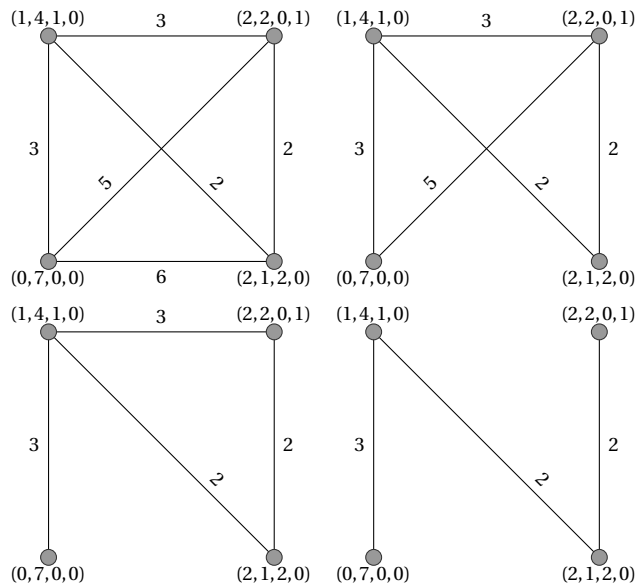
\begin{center}
\begin{tikzpicture}[y=.3cm, x=.3cm,font=\scriptsize]
\draw (0,0) -- (10,10);
\draw (0,0) -- (0,10);
\draw (0,0) -- (10,0);
\draw (10,0) -- (10,10);
\draw (10,0) -- (0,10);
\draw (0,10) -- (10,10);

\filldraw[fill=black!40,draw=black!80] (0,0) circle (3pt)    node[anchor=north] {$(0,7,0,0)$};

\filldraw[fill=black!40,draw=black!80] (0,10) circle (3pt)    node[anchor=south] {$( 1, 4, 1, 0 )$};

\filldraw[fill=black!40,draw=black!80] (10,0) circle (3pt)    node[anchor=north] {$( 2, 1, 2, 0 )$};

\filldraw[fill=black!40,draw=black!80] (10,10) circle (3pt)    node[anchor=south] {$(2, 2, 0, 1)$};

\node [above] at (5,10) {$3$};

\node [below] at (5,0) {$6$};

\node [right] at (10,5) {$2$};

\node [left] at (0,5) {$3$};

\node [above, rotate=45] at (3,3) {$5$};

\node [above, rotate=-45] at (7,3) {$2$};

\end{tikzpicture}
\begin{tikzpicture}[y=.3cm, x=.3cm,font=\scriptsize]
\draw (0,0) -- (10,10);
\draw (0,0) -- (0,10);
\draw (10,0) -- (10,10);
\draw (10,0) -- (0,10);
\draw (0,10) -- (10,10);

\filldraw[fill=black!40,draw=black!80] (0,0) circle (3pt)    node[anchor=north] {$(0,7,0,0)$};

\filldraw[fill=black!40,draw=black!80] (0,10) circle (3pt)    node[anchor=south] {$( 1, 4, 1, 0 )$};

\filldraw[fill=black!40,draw=black!80] (10,0) circle (3pt)    node[anchor=north] {$( 2, 1, 2, 0 )$};

\filldraw[fill=black!40,draw=black!80] (10,10) circle (3pt)    node[anchor=south] {$(2, 2, 0, 1)$};

\node [above] at (5,10) {$3$};


\node [right] at (10,5) {$2$};

\node [left] at (0,5) {$3$};

\node [above, rotate=45] at (3,3) {$5$};

\node [above, rotate=-45] at (7,3) {$2$};

\end{tikzpicture}

\begin{tikzpicture}[y=.3cm, x=.3cm,font=\scriptsize]
\draw (0,0) -- (0,10);
\draw (10,0) -- (10,10);
\draw (10,0) -- (0,10);
\draw (0,10) -- (10,10);

\filldraw[fill=black!40,draw=black!80] (0,0) circle (3pt)    node[anchor=north] {$(0,7,0,0)$};

\filldraw[fill=black!40,draw=black!80] (0,10) circle (3pt)    node[anchor=south] {$( 1, 4, 1, 0 )$};

\filldraw[fill=black!40,draw=black!80] (10,0) circle (3pt)    node[anchor=north] {$( 2, 1, 2, 0 )$};

\filldraw[fill=black!40,draw=black!80] (10,10) circle (3pt)    node[anchor=south] {$(2, 2, 0, 1)$};

\node [above] at (5,10) {$3$};


\node [right] at (10,5) {$2$};

\node [left] at (0,5) {$3$};


\node [above, rotate=-45] at (7,3) {$2$};

\end{tikzpicture}
\begin{tikzpicture}[y=.3cm, x=.3cm,font=\scriptsize]
\draw (0,0) -- (0,10);
\draw (10,0) -- (10,10);
\draw (10,0) -- (0,10);

\filldraw[fill=black!40,draw=black!80] (0,0) circle (3pt)    node[anchor=north] {$(0,7,0,0)$};

\filldraw[fill=black!40,draw=black!80] (0,10) circle (3pt)    node[anchor=south] {$( 1, 4, 1, 0 )$};

\filldraw[fill=black!40,draw=black!80] (10,0) circle (3pt)    node[anchor=north] {$( 2, 1, 2, 0 )$};

\filldraw[fill=black!40,draw=black!80] (10,10) circle (3pt)    node[anchor=south] {$(2, 2, 0, 1)$};



\node [right] at (10,5) {$2$};

\node [left] at (0,5) {$3$};


\node [above, rotate=-45] at (7,3) {$2$};
\end{tikzpicture}
\end{center}
Thus the catenary degree of $77$ is $3$.
\end{example}

Observe that in Remark \ref{remark-pres}, we obtained chains joining any two factorizations of the same element, just using translations of elements in a presentation. Since distances are not translation-sensitive, we only have to care on how to go from the first component to the second in a relation in a presentation. It follows (see \cite{cat-tame}) that 
 \[
\mathsf c(M)=\max\big\{ \mathsf c(b)\mid b\in \mathrm{Betti}(M)\big\}.
\]
This gives a computational procedure to compute the catenary degree of any affine semigroup $M$.

\begin{example}
Let us recover Example \ref{pres-min}, $M=\langle (2,0),(0,2),(1,1),(2,1)\rangle$. We already know that $\mathrm{Betti}(M)=\{(2,2),(2,4)\}$. 
\begin{verbatim}
gap> a:=AffineSemigroup([2,0],[0,2],[1,1],[1,2]);;
gap> gens:=GeneratorsOfAffineSemigroup(a);                     
[ [ 0, 2 ], [ 1, 1 ], [ 1, 2 ], [ 2, 0 ] ]
gap> betti:=BettiElementsOfAffineSemigroup(a);
[ [ 2, 2 ], [ 2, 4 ] ]
gap> List(betti,b->FactorizationsVectorWRTList(b,gens));
[ [ [ 1, 0, 0, 1 ], [ 0, 2, 0, 0 ] ], 
  [ [ 2, 0, 0, 1 ], [ 1, 2, 0, 0 ], [ 0, 0, 2, 0 ] ] ]
gap> List(last,CatenaryDegreeOfSetOfFactorizations);
[ 2, 3 ]
gap> CatenaryDegreeOfAffineSemigroup(a);
3
\end{verbatim}
\end{example}

So far we do not know of a procedure to compute the (finite) set $\{\mathsf c(m)\mid m\in M\}$. It is known that for numerical semigroups, the catenary degree is also eventually periodic, but unfortunately no bounds for this periodicity are known (\cite{ct-per}). For half-factorial monoids it can be shown (see \cite[Theorem 2.3]{hom}) that 
\[
\{\mathsf c(m)\mid m\in M\}=\{\mathsf c(m)\mid m\in \mathrm{Betti}(M)\}.
\]

For numerical semigroups, in light of Section \ref{sec:apery} (see also \cite[Corollary 3]{cat-tame-ns}), 
\[
\mathsf c(M)=\max\big\{\mathsf c(m)\mid m\in \{m_2,\ldots,m_e\}+(\mathrm{Ap}(M,m_1)\setminus\{0\}) \big\},
\]
and so in this setting it is not needed to compute $\mathrm{Betti}(M)$.

\subsection{Monotone, equal and homogeneous catenary degrees} We can obtain different flavored catenary degrees if we impose conditions on the definition of $N$-chain. For instance if we enforce the chain of factorizations to have nondecreasing lengths we obtain the definition of \emph{monotone catenary degree}. 

We can also ask the lengths to be all equal, and then we have \emph{equal catenary degree}. 

Finally we can also impose that the lengths in the chain are not larger than the maximum of the lengths of the ends of the chain, obtaining in this way the \emph{homogeneous catenary degree}. 

Let $M$ be an affine semigroup, $M\subseteq \mathbb N^k$. Let $m\in M$, we write $(m,1)\in \mathbb N^{k+1}$ for the element with the first coordinates the coordinates of $m$ and last coordinate equal to 1 (we have appended a 1 at the ``end'' of $m$). Assume that $M$ is minimally generated by $\{m_1,\ldots, m_e\}$. Set 
\[
M^{eq}=\langle (m_1,1),\ldots, (m_e,1)\rangle
\]
and 
\[
M^{hom}=\langle (m_1,1),\ldots, (m_e,1),(0,1)\rangle.
\]
Both $M^{eq}$ and $M^{hom}$ are half-factorial monoids (see \cite{hom}).
It is easy to prove that the equal catenary degre of $M$ corresponds with the catenary degree of $M^{eq}$, and that the homogeneous catenary degree of $M$ is precisely $M^{hom}$.  This provides a way to compute both homogeneous and equal catenary degrees.  

In order to compute the monotone catenary degre of $M$, it can be derived from \cite{phi-2} that we have to look at the projections in the first $k$ coordinates of the primitive elements of $M^{hom}$ (see \cite[Chapter 3]{tesis-alfredo}), and then take the maximum of the monotone catenary degrees of these elements. The monotone catenary degree of $m$ is the maximum of the equal and adjacent catenary degree of $m$, where the {adjacent catenary degree} of $m$ is defined as follows: let $\mathsf L(m)=\{l_1<\cdots <l_{r}\}$, and for every $i\in\{1,\ldots,r\}$ denote by $\mathsf Z_{l_i}(m)$ the set of factorizations of $m$ with length $l_i$; the \emph{adjacent catenary degree} of $m$ is the maximum of the distances $\mathrm d(Z_{l_i},Z_{l_{i+1}})$, $i\in\{1,\ldots,r-1\}$.

\begin{example}
Let us use \texttt{numericalsgps} to compute the catenary degrees of $\langle 10,17,24,31,43\rangle$.

\begin{verbatim}
gap> s:=NumericalSemigroup(10,17,24,31,43);
<Numerical semigroup with 5 generators>
gap> MinimalGeneratingSystem(s);           
[ 10, 17, 24, 31, 43 ]
gap> CatenaryDegreeOfNumericalSemigroup(s);
6
gap> HomogeneousCatenaryDegreeOfNumericalSemigroup(s);
11
gap> MonotoneCatenaryDegreeOfNumericalSemigroup(s);
11
gap> EqualCatenaryDegreeOfNumericalSemigroup(s);
11
\end{verbatim}
\end{example}

\subsection{Tame degree}
Assume that $M$ is an affine semigroup generated by $\{m_1,\ldots, m_e\}$, and let $m$ in $M$ and $x\in \mathsf Z(m)$. If there exists $n_1,\ldots, n_e\in \mathbb N$ such that $m-\left(\sum_{i=1}^em_i\right)\in M$, then there must be $y=(y_1,\ldots,y_e)\in \mathsf Z(m)$ such that $y-(n_1,\ldots,n_e)\in \mathbb N^e$. We want to know the smallest possible distance at which we can find such an $y$. This is the idea of tame degree. We are mostly interested in the case $\sum_{i=1}^en_i m_i=m_j$ for some $j\in\{1,\ldots, e\}$.

The \emph{tame degree} of $m$ with respect to $m_i$, $\mathsf t(m, m_i)$, is the least nonnegative integer $t$ such that for every $z\in\mathsf Z(m)$, there exists $z'\in\mathsf Z(m)$ with $z'-\mathbf e_i\in \mathbb N^k$ (or in other words, the $i$th coordinate of $z'$ is nonzero) and $\mathrm d(z,z')\le t$. The tame degree of $M$ with respect to $m_i$, $\mathsf t(M,m_i)$, is the supremum (maximum in this setting, \cite{cat-tame}) of all the tame degrees of the elements of $m_i+M$ with respect to $m_i$. 

The tame degree of $M$, $\mathsf t(M)$, is the maximum of the tame degrees of $S$ with respect to all the atoms (affine semigroups are tame and locally tame, \cite{g-hk}). The tame degree of $M$ can be computed by means of the tame degrees of the primitive elements of $M$ (\cite{cat-tame}). Recently, a faster approach has been described in \cite{tame}. Set $\mathcal M_i=\mathrm{Minimals}_\le \mathsf Z(m_i+M)$ and $M_i=\{\varphi(z)\mid z\in \mathcal M_i\}$. By Dickson's lemma, $\mathcal M_i$ and $M_i$ have finitely many elements. Moreover,
\[ 
\mathsf t(M,m_i)=\max\big\{ \mathsf t(m,m_i)\mid m\in M_i\big\}.
\]
In \cite{minimales} there is a procedure to compute $M_i$ (indeed the set of expressions of any ideal of $M$, not just principal ideals). By using \cite{normaliz} or \cite{4ti2} (or any integer linear programming package) we can also compute this directly as in the following example.

\begin{example}\label{ideal-principal}
Let us compute the set $M_1$ for $M=\langle (2,0),(0,2),(1,1),(1,2)\rangle$. We need to find the expressions in $(2,0)+M$. This corresponds with the $(x,y,z,t)\in \mathbb N^4$ such that 
\[
\begin{pmatrix}
2 & 0 & 1 & 1 \\
0 & 2 & 1 & 2
\end{pmatrix}
\begin{pmatrix}
x \\ y \\ z \\ t
\end{pmatrix}
=
(2,0)+
\begin{pmatrix}
2 & 0 & 1 & 1 \\
0 & 2 & 1 & 2
\end{pmatrix}
\begin{pmatrix}
x' \\ y' \\ z' \\ t'
\end{pmatrix}
\]
for some $(x',y',z',t')\in \mathbb N^4$. And this is a system of two equations and eight unknowns.	We use the package \texttt{4ti2gap} to solve this.

\begin{verbatim}
gap> m:=[[2,0,1,1,-2,0,-1,-1],[0,2,1,2,0,-2,-1,-2]];             
[ [ 2, 0, 1, 1, -2, 0, -1, -1 ], [ 0, 2, 1, 2, 0, -2, -1, -2 ] ]
gap> problem:=["mat",m,"rhs",[[2,0]],"sign",[[1,1,1,1,1,1,1,1]]];
[ "mat", [ [ 2, 0, 1, 1, -2, 0, -1, -1 ], [ 0, 2, 1, 2, 0, -2, -1, -2 ] ], 
"rhs", [ [ 2, 0 ] ], "sign", [ [ 1, 1, 1, 1, 1, 1, 1, 1 ] ] ]
gap> ZSolve4ti2(problem);
rec( zhom := [ [ 1, 0, 0, 2, 0, 0, 4, 0 ], [ 0, 0, 4, 0, 1, 0, 0, 2 ], 
[ 1, 2, 0, 0, 0, 0, 0, 2 ], [ 0, 0, 1, 0, 0, 0, 1, 0 ], 
[ 0, 0, 0, 1, 0, 0, 0, 1 ], [ 0, 1, 2, 0, 0, 0, 0, 2 ], 
[ 0, 0, 2, 0, 1, 1, 0, 0 ], [ 0, 1, 0, 0, 0, 1, 0, 0 ], 
[ 0, 0, 0, 2, 1, 2, 0, 0 ], [ 1, 0, 0, 0, 1, 0, 0, 0 ], 
[ 1, 1, 0, 0, 0, 0, 2, 0 ], [ 0, 0, 0, 2, 0, 1, 2, 0 ] ], 
zinhom := [ [ 0, 0, 4, 0, 0, 0, 0, 2 ], [ 0, 0, 2, 0, 0, 1, 0, 0 ], 
[ 1, 0, 0, 0, 0, 0, 0, 0 ], [ 0, 0, 0, 2, 0, 2, 0, 0 ] ] )
\end{verbatim}
This in particular means that 
\begin{multline*}
	\mathsf Z((2,0)+M)=\{(0,0,4,0), (0,0,2,0), (1,0,0,0), (0,0,0,2)\}\\+\langle (1,0,0,1), (0,0,4,0), (1,2,0,0), 
(0,0,1,0), (0,0,2,0), (0,1,0,0), (0,0,0,2),(1,0,0,0), (1,1,0,0)\rangle
\end{multline*}
And thus $\mathrm{Minimals}_\le \mathsf Z((2,0)+M)=\{(0,0,2,0),(1,0,0,0),(0,0,0,2)\}$. Hence $M_1=\{ (2,2), (2,0), (2,4)\}$.
\end{example}

If $M$ is a full affine semigroup (for instance in the case of block monoids), then the elements in $\mathcal M_i$ can be computed using \cite[Corollary 3.5]{b-gs-g}. In this case $M_i$ is the set of minimal nonnegative integer solutions of 
\[
(m_1\mid \cdots \mid m_e)x^T \ge m_i.
\]

\begin{example}
Let us compute as explained in \cite{tame} the tame degree of $\mathcal B\left(\mathbb Z_2^3\right)$. 	
\begin{verbatim}
gap> c:=[ [ 0, 0, 1 ], [ 0, 1, 0 ], [ 0, 1, 1 ], [ 1, 0, 0 ], [ 1, 0, 1 ], 
[ 1, 1, 0 ], [ 1, 1, 1 ] ];
gap> a:=AffineSemigroup("equations",[TransposedMat(m),[2,2,2]]);;
gap> TameDegreeOfAffineSemigroup(a);
4
\end{verbatim}
\end{example}

For numerical semigroups we have a similar behavior as in the catenary degree. The tame degree is reached in an element that has to do with Ap\'ery sets (\cite[Theorem 16]{cat-tame-ns})):
\[
\mathsf t(M)=\max\left\{\mathsf t(m) ~\Big|~ m\in \{m_1,\ldots,m_e\}+\left(\bigcup_{i=1}^e \mathrm{Ap}(M,m_i)\setminus\{0\}\right)\right\}.
\]
For small generators, the above formula is faster than computing minimal factorizations in principal ideals (or if we do not have software to solve linear Diophantine equations over the set of nonnegative integers at hand).

\section{$\omega$-primality}

Let $M$ be an affine semigroup. Define on $M$ the following binary relation: $m\le_M m'$ if $m'-m\in M$. This relation is an order relation (the translation of divisibility to additive notation). We say that $m\in M$ is \emph{prime} if whenever $m\le_M m'+m''$ for some $m',m''\in M$, either $m\le_M m'$ or $m\le_M m''$. Any prime element must be an atom. But it may happen that no atom is prime (this holds in any nontrivial numerical semigroup). The $\omega$-primality is meant to determine how far is an element from being prime.

The \emph{$\omega$-primality} of $m$ in $M$, denoted $\omega(m)$, is the least positive integer $N$ such that whenever $m\le_M a_1+\cdots+a_n$ for some $a_1,\ldots, a_n\in M$, then $m\le_M a_{i_1}+\cdots+a_{i_N}$ for some  $\{i_1,\ldots,i_N\}\subseteq \{1,\ldots, n\}$. 

According to this definition an element is prime provided that its $\omega$-primality is one.

Notice that by definition, $m\le_M m'$ if and only if $m'$ is in the principal ideal $m+M$.  Hence principal ideals play a fundamental role in the computation of $\omega$-primality (as in the calculation of the tame degree). Indeed in \cite[Proposition 3.3]{b-gs-g} it is shown that 
\[ \omega(m)= \max\big\{\lvert x\rvert  ~\big|~  x\in \mathrm{Minimals}_{\leq}(\mathsf Z(m+M))\big\}.\]
 In \cite{omega-comp} the above formula together with the algorithm presented in \cite{minimales} is used to compute the $\omega$-primality of an element in an affine semigroup. One can also proceed as in Example \ref{ideal-principal} and use for instance \texttt{Normaliz} or \texttt{4ti2}.

The omega primality of $M$, if $M$ is minimally generated by $\{m_1,\ldots,m_e\}$, is defined as $\omega(M)$ as the maximum of $\{\omega(m_1),\ldots, \omega(m_e)\}$. Note that $\omega(m)$ with $m$ running in $M$ is not upper bounded in general.

\begin{example}
According to Example \ref{ideal-principal}, for $M=\langle (2,0),(0,2),(1,1),(1,2)\rangle$, we have $\omega((2,0))=2$. Let us double check it with the \texttt{numericalsgps} package.
\begin{verbatim}
gap> a:=AffineSemigroup([2,0],[0,2],[1,1],[1,2]);;
gap> OmegaPrimalityOfElementInAffineSemigroup([2,0],a);
2
gap> OmegaPrimalityOfAffineSemigroup(a);
4
\end{verbatim}
\end{example}

For numerical semigroups we obtain a similar construction as for the tame degree (as expected, since we are using roughly the same elements in the calculations). In \cite[Remarks 5.9]{b-gs-g} it is shown that if we are looking for minimal factorizations in $\mathsf Z(m+M)$, then we only have to search for factorizations of the elements of the form $m+w$ with $w\in \mathrm{Ap}(M,m_i)$ for some $i\in\{1,\ldots, e\}$. In \cite{b-on-p} an improved method that also uses Ap\'ery sets is given (its actually the one that uses the package \texttt{numericalsgps}; see \texttt{contributions.gi} in the package).

\begin{example}
Let us compare the timings for $S=\langle 10,17,24,31,43\rangle$.
\begin{verbatim}
gap> s:=NumericalSemigroup(10,17,24,31,43);;
gap> OmegaPrimalityOfNumericalSemigroup(s);time;
11
13
gap> a:=AsAffineSemigroup(s);;                     
gap> OmegaPrimalityOfAffineSemigroup(a);time;    
11
3654
\end{verbatim}
(The timings are in milliseconds.)

If the generators are larger, then the principal ideal approach is better.
\begin{verbatim}
gap> s:=NumericalSemigroup(201,223,357);;
gap> OmegaPrimalityOfNumericalSemigroup(s);time;
75
32245
gap> a:=AsAffineSemigroup(s);;
gap> OmegaPrimalityOfAffineSemigroup(a);time;   
75
1934
\end{verbatim}
\end{example}

\end{document}